\documentclass[11pt,a4paper,leqno]{amsart}

\usepackage{fullpage,amsmath,amsthm,amsfonts,amscd,
latexsym,amssymb,eucal,curves}
  \usepackage[all]{xy}
\usepackage[dvips]{graphicx}
\usepackage{epsfig,subfigure}
\usepackage{geometry,psfrag}
\newcommand{\ZZ}{\mathbb{Z}}

\newcommand{\CC}{\mathbb{C}}

\newcommand{\PP}{\mathbb{P}}
\newcommand{\NN}{\mathbb{N}}

\newcommand{\HH}{\mathbb{H}}
\newcommand{\GG}{\mathbb{G}}

\newcommand{\LL}{\mathbb{L}}

\newcommand{\TT}{\mathbb{T}}
\newcommand{\MM}{\mathbb{M}}
\newcommand{\UU}{\mathbb{U}}

\newcommand{\QQ}{\mathbb{Q}}
\newcommand{\RR}{\mathbb{R}}

\newcommand{\VV}{\mathbb{V}}

\newcommand{\Jac}{{\rm Jac} \,}

\newcommand{\Aff}{{\rm Aff}}

\newcommand{\cLL}{{\mathcal L}}
\newcommand{\cOO}{{\mathcal O}}
\newcommand{\cEE}{{\mathcal E}}
\newcommand{\cFF}{{\mathcal F}}

\newcommand{\cYY}{{\mathcal Y}}
\newcommand{\cXX}{{\mathcal X}}

\newcommand{\SL}{{\rm SL}}

\newcommand{\diag}{{\rm diag}}

\newcommand{\opt}{{\rm opt}}

\newcommand{\sms}{\smallsetminus}

\DeclareMathOperator{\Sp}{Sp}

\newcommand{\ol}{\overline}

\newcommand{\wtx}{\widetilde{\cXX}}

\newcommand{\Hg}{{\rm Hg}}

\newcommand{\TM}{Teich\-m\"ul\-ler\ }

\newtheorem{Defi}{Definition}[section]

\newtheorem{Rem}[Defi]{Remark}

\newtheorem{Prop}[Defi]{Proposition}
\newtheorem{Lemma}[Defi]{Lemma}
\newtheorem{Cor}[Defi]{Corollary}
\newtheorem{Thm}[Defi]{Theorem}

\makeatletter
\renewcommand{\subsection}{\@startsection{subsection}{2}%
        {\z@}{-1.25ex plus -1ex minus-.2ex}{-1em}{\bf}}
\makeatother

\begin{document}{\large}
\title{Shimura- and Teichm\"uller curves}
\author{Martin M\"oller}

\begin{abstract}
We classify curves in the moduli space of curves $M_g$
that are both Shimura- and Teichm\"uller curves: 
For both $g=3$ and $g=4$ there exists precisely one such
curve, for $g=2$ and $g \geq 6$ there are no such curves.
\newline
We start with a Hodge-theoretic description of 
Shimura curves and of Teichm\"uller curves that reveals similarities
and differences of the two classes of curves.  
The proof of the classification relies on the geometry of 
square-tiled coverings and on estimating
the numerical invariants of these particular fibered surfaces.
\newline
Finally we translate our main result into a classification
of Teichm\"uller curves with totally degenerate Lyapunov
spectrum.
\date{\today}
\end{abstract}
\thanks{Supported by the 
DFG-Schwerpunkt ``Komplexe Mannigfaltigkeiten'' and by the MPIM, Bonn}
\maketitle

\section*{Introduction}

A Teichm\"uller curve is an algebraic curve $C \to M_g$ in the moduli 
space of curves, which is the image of a holomorphic geodesic for the
Teichm\"uller (equivalently: Kobaya\-shi-) metric $\tilde{\j}: \HH \to T_g$ 
in Teichm\"uller space. Only rarely geodesics in Teichm\"uller space 
map to algebraic curves in $M_g$. 
\par
By Teichm\"uller's theorems a geodesic $\tilde{\j}$ is generated
by a Riemann surface $X$ with a Teichm\"uller marking together
with a quadratic differential $q$ on $X$. By \cite{Kr81} (or
\cite{McM03} Theorem.\ 4.1) the composition $C \to M_g \to A_g$ is a
geodesic for the Kobayashi metric precisely if $q = \omega^2$
for some holomorphic $1$-form $\omega$ on $X$. We deal here
exclusively with these Teichm\"uller curves.
\par
A Shimura curve of Hodge type is a curve $\HH/\Gamma =:C \to A_g$ 
in the moduli space of abelian varieties of dimension $g$ that
is totally geodesic for the Hodge metric on $A_g$ and that contains
a CM point. Such curves are automatically algebraic. See 
Section~\ref{Shimura} for a group-theoretic definition of Shimura curves 
\par
Both Shimura curves and Teichm\"uller curves can be characterized by
their variation of Hodge structures (VHS). See Section~\ref{Shimura} and 
Section~\ref{Teich} for similarities and differences between the two 
sorts of curves. 
One also could take the Theorems~\ref{VHSShimura} and \ref{VHSTeich} 
stated there as definition of Shimura- resp.\ Teichm\"uller curves.
\par
Here we investigate whether there are curves with {\em both} 
properties, i.e.\ 
curves in $M_g$ that are Teichm\"uller curves and, when considered
in $A_g$, Shimura curves. We call them  { \em ST-curves} for short.
The motivation for studying them is threefold:
\newline
First, the VHS of a Teichm\"uller curve consists of a sub-local system $\LL$
that is maximal Higgs (see Section~\ref{Teich} for the definition), 
its conjugates and some rest $\MM$. Not much is known about $\MM$,
except that $\MM$ is a 'contraction' of $\LL$ (see \cite{McM03})
in the following sense: If $\gamma \in \pi_1(C)$ is hyperbolic,
the largest eigenvalue of $\gamma$ acting on a fiber of $\LL$ is strictly
bigger than the the largest eigenvalue of $\gamma$ acting on a 
fiber of $\MM$. Shimura curves correspond to the case where
$\MM$ is as 'small' as possible, namely trivial.
\newline
Second, they are curves with 'few' singular fibers: Beauville
and Tan have shown that semistable fibrations over $\PP^1$
have at least $4$ singular fibers resp.\ $5$, 
if the fiber genus is at least two. 
In \cite{TaTuZa04} the authors study semistable fibrations
over $\PP^1$ with $5$ and $6$ singular fibers. Instead of sticking
to the basis $\PP^1$ one could admit any base and ask for the
following: Classify fibrations that have few singular fibers with 
proper Jacobian compared to the number of singular fibers with
non-smooth fibrations! The curves that are both Shimura- and
Teichm\"uller- have {\em no} singular fibers with 
smooth Jacobian by Proposition~\ref{noJsmoothsing}. 
\newline
The third motivation concerns the Lyapunov exponents of the
of lift of the \TM geodesic flow to the Hodge bundle
(see Section~\ref{sec:Lyap}).
For the natural measure supported on whole strata of the Hodge bundle
the spectrum of Lyapunov exponents has recently been
shown to be simple (\cite{Fo02}, \cite{AV07}). In contrast, the
Lyapunov spectrum for the natural measure supported on an 
individual Teichm\"uller curves can
be maximally degenerate, i.e.\ all but the top and bottom
Lyapunov exponent are zero. We will show in Section~\ref{sec:Lyap}
that these curves are precisely the ST-curves.
\par
As explained in Section~\ref{Teich} we abuse the notion 
Teichm\"uller curve for \'etale coverings of what we have defined
above. Our main result is:
\par
{\bf Theorem~\ref{Main}:} {\em For $g=2$ and $g \geq 6$ there are
no ST-curves.
\newline
In both $M_3$ and in $M_4$ there is only one 
ST-curve. Its 
universal family is given by
$$ y^4 = x(x-1)(x-t) $$
in $M_3$ and respectively in $M_4$ by
$$ y^6 = x(x-1)(x-t). $$
In both cases  $t \in \PP^1 \sms \{0,1,\infty\}$. }
\par
We remark that the combinatorial discussion in the cases $g=4$
and $g=5$ was incomplete when earlier versions of this
paper circulated. The example in $M_4$ was shown to
the author by Forni and Matheus (\cite{FM08}).
There are many constraints for the existence of a  ST-curve
in genus $5$, see Corollary~\ref{NoeCor} and Corollary~\ref{eqlenalldir}. We
conjecture that there is no ST-curve in $M_5$.
\par
The classification of ST-curves
reminds of the conjectures of Colemen and Andr\'e-Oort
that together imply that there should be no Shimura curve
{\em generically} lying in the moduli space of curves $M_g$ 
for large $g$. The main theorem can be rephrased saying that there is
no non-compact Shimura curve $C$ lying {\em entirely} in $M_g$ 
for $g \geq 4$. We emphasize that 'entirely' refers not only
to the non-existence of singular fibers with proper Jacobian
but also that the Torelli map is unramified if $C$ should happen
to contain points in the image of the hyperelliptic locus. 
After a first version of this paper appeared, it was shown 
in \cite{ViZu06} that the hypothesis 'non-compact' can be
removed. Together these results imply:
\par
{\bf Corollary:} There is no Shimura curve in $A_g$ for
$g \geq 6$ that lies entirely in $M_g$.
\par
Unfortunately, for Shimura curves lying only generically in the 
moduli space of curves $M_g$ none of our  techniques
of flat geometry apply, since such a curve can never be
a Teichm\"uller curve.
\par
Using the characterization in terms of VHS one deduces that a 
Teichm\"uller curve
can be Shimura only if it arises as square-tiled
covering, i.e.\ with affine group $\Gamma$ commensurable to $\SL_2(\ZZ)$. 
The notions are explained in Section~\ref{Teich} and this key step
in stated as Corollary~\ref{STimpliesr1}.
Given any element of the affine group of a square-tiled covering, one can 
effectively compute the action of the fundamental group on  $H^1(X_0,\ZZ)$ of a fiber $X_0$ 
and hence decide whether the covering gives a Shimura curve or not. 
But since it is hardly
clear how the geometry of the covering translates into properties of the
action of the affine group, it seems not tractable to identify Shimura curves
in this way.
\par
Our proof is thus rather indirect: Being both a Shimura and a 
Teichm\"uller curve imposes strong conditions on the geometry of 
the degenerate fibers.
Exploiting these plus some geometric considerations on
flat surfaces suffices to treat the case $g \leq 3$. The non-existence 
for larger $g$ follows from calculating both sides of the Noether 
formula for the fibered surface $f:\cXX \to C$ using that we can
work in fact of a modular curve $C=X(d)$. We end up with
a finite list of possible covering degrees and types of
zeros for $\omega$ in genus $g=4$ and $g=5$. To show that
these cases in fact do not exist for $g=4$, we translate the geometric
information about degenerate fibers into conditions on the
monodromy representation of the square-tiled coverings.
\par
The Shimura- and Teichm\"uller curve in $M_3$ above appears in various 
guises in the literature: The author heard about its properties
independently from Herrlich and Schmitth\"usen (\cite{HeSc04}) and
Forni (\cite{Fo06}).
If we consider it over the modular curve $X(4) \cong \PP^1 \sms \{6 \; {\rm
points} \}$ instead of over $X(2)$ then it admits a stable model.
The fibration hence fits in the context of \cite{TaTuZa04}, where
it is shown that the total space over $X(4)$ is a K3-surface.
\newline
The decomposition of the Jacobian of $y^N=x(x-1)(x-t)$ has been studied
by many authors, e.g.\ \cite{Wo88}.
Standard criteria for the monodromy groups of hypergeometric differential
equation  imply that for $N=4$ and $N=6$ the family defines indeed 
a Shimura curve. 
\newline
In \cite{Gu03} the family for $N=4$ is studied from the Arakelov viewpoint
and \cite{HeSc04} analyze how this Teichm\"uller curve intersects
other Teichm\"uller curves.
\par
The author thanks Frank Herrlich, Gabi Schmitth\"usen 
and Eckart Viehweg for many fruitful discussions.
He also thanks Giovanni Forni and Carlos Matheus Santos
for sharing their insights on the cyclic coverings. 
Thomas Fischbacher helped the author with computer search of examples
in genus 5, using \cite{Fi02}. Finally, the author thanks the
referees for their comments that helped to improve the exposition
of the result.

\section{Shimura curves} \label{Shimura}

In this section we give a characterization of Shimura curves
by their variation of Hodge structures. The hard part of this
characterization, Arakelov equality implies Shimura curve, is
the content of \cite{ViZu04}. These authors were aware of the
converse implication, which we need later on, but a proof does
not appear in \cite{ViZu04}, nor to the author's knowledge
anywhere else in the literature. We thus present this proof assuming
that the reader is familiar with the basic notions of algebraic groups. 
\par
For a reader with background in Teichm\"uller we suggest to only compare
the results of Theorem~\ref{VHSShimura} and Theorem~\ref{VHSTeich}. 
To make this possible in a self-contained way, we recall some 
generalities on variations of Hodge structures. Let
 $f: A \to C$ be a family of abelian varieties over a smooth curve
$C$, let $\ol{C}$ be its smooth compactification and $S := \ol{C} \sms C$.
We denote by $\VV = R^1 f_* \ZZ$ the corresponding local system. 
Take any completion $f: \ol{A} \to \ol{C}$ with  $\ol{A}$ smooth. 
Then the Deligne extension $\VV \otimes \cOO_C$ to $\ol{C}$ carries
a Hodge filtration 
$$ 0 \subset f_* \omega_{\ol{A}/\ol{C}} \subset 
(\VV \otimes \cOO_C)_{\rm ext}.$$ 
The local system $\VV$ and the Hodge filtration form a variation of Hodge
structures (VHS) of weight one. 
The graded pieces of the Hodge filtration together with 
the quotient map induced by Gauss-Manin connection form a Higgs-bundle
$$(\cEE^{1,0} \oplus \cEE^{0,1}, \Theta:\cEE^{1,0} \to \cEE^{0,1} 
\otimes \Omega_{\ol{C}}(\log S) ), $$
where $\cEE^{1,0}= f_* \omega_{\ol{A}/\ol{C}}$ and 
$\cEE^{0,1} = R^1 f_* \cOO_{\ol{A}}$.
\par
Given such a Higgs bundle, the subbundle with Higgs field 
$\Theta=0$ can be split off as a direct summand (see \cite{Ko87}).   
The remaining part $(\cFF^{1,0} \oplus \cFF^{0,1})$ satisfies the
Arakelov inequality 
$$ 2\cdot \deg(\cFF^{1,0}) \leq {\rm rank}(\cFF^{1,0})\cdot 
\deg(\Omega_{\ol{C}}^1(\log(S))). $$
\par
A VHS $\LL$ of rank two is called {\em maximal Higgs} (\cite{ViZu04}), 
if the corresponding Higgs bundle satisfies Arakelov equality or
equivalently, if $\Theta$ is an isomorphism.
\par
\begin{Defi}
A {\em Shimura datum} consists of 
\begin{itemize}
\item{i)} a reductive algebraic group $G$ defined over $\QQ$,
\item{ii)} a vector space $V$ with lattice $L$ and a symplectic
paring $Q$, which is integral on $L$,
\item{iii)} a faithful representation $\rho: G \to {\rm Sp}(V,Q)$,
\item{iv)} a complex structure $\varphi_0: S^1 = \{z \in \CC: |z| =1\} \to
{\rm Sp}(V,Q)_\RR$ such that $Q(x,\varphi_0(i)x)>0$ for all $x \neq 0$,
\end{itemize}
such that
$$ \rho(G) \;\text{is normalized by} \;\varphi_0(S^1) \eqno(H_1)$$
\par
Given these data, $K_\RR =\{g \in G_\RR | 
\rho(g)\varphi_0= \varphi_0 \rho(g)\}$ is a maximal
compact subgroup of $G_\RR$. The quotient ${\mathcal X} =
G_\RR/K_\RR$ is a bounded
symmetric domain and $V_\RR/ L$ with complex structure 
$\rho(g)\varphi_0\rho(g)^{-1}$
defines a holomorphic family of abelian varieties over $D$.
\par
A {\em Shimura variety}
is the inverse system $\{{\mathcal X}/\Gamma\}_\Gamma$ where $\Gamma$ runs
over the torsion-free congruence subgroups of $G(\QQ)$ such
that $\rho(\Gamma)$ preserves $L$. 
\par
The {\em Hodge group} $\Hg := \Hg(\varphi_0)$ is the smallest $\QQ$-algebraic
subgroup of ${\rm Sp}(V,Q)$ that contains the image $\varphi_0(S^1)$.
\par
The Shimura datum (and the Shimura variety) is {\em of Hodge type} 
if $G = \Hg(\varphi_0)$.
\end{Defi}
\par
\par
For a sufficiently small arithmetic subgroup 
$\widetilde{\Gamma}$ of ${\rm Sp}_{2g}$ and $\Gamma= \widetilde{\Gamma} 
\cap \rho(G)$
we have a map
$$j: {\mathcal X}/\Gamma\ =:C \to A_g^{\widetilde{\Gamma}}$$
into the moduli space
of abelian varieties with some level structure $\widetilde{\Gamma}$.
\par
From an embedding $j$ one obviously regains the whole Shimura
datum attached to $C$ and $A_g$ and we will henceforth call
$j$ a (representative of a) Shimura variety. From now on we
exclusively deal with {\em Shimura curves}, i.e.\ 
Shimura varieties of dimension one. In this case the symmetric
domain ${\mathcal X}$ will be the upper half plane $\HH$.
\par
Let $f: A \to C$ be the pullback of the universal family
over $A_g^{\widetilde{\Gamma}}$ to $C$. We will always suppose that
the monodromies around the cusps of Shimura- and Teichm\"uller
curves are unipotent, replacing $\Gamma$ by a subgroup of finite
index if necessary. 
\par

Let $X_{\Hg}$ be the $\Hg(\RR)^+$-conjugacy class in $
{\rm Hom}_{{\rm alg. grp}/\RR}(  {\rm Res}_{\CC/\RR}\GG_m , \Hg_\RR)
$
containing $\varphi \circ ({\rm Res}_{\CC/\RR}\GG_m \to S^1)$. Here $+$
denotes the topological connected component. 
In the definition of a Shimura datum we closely followed Mumford (\cite{Mu69},
see also \cite{Sa80} for the condition $(H_1)$). 
For the reader's convenience we show for comparison that $(\Hg, X_\Hg)$ 
is a Shimura datum in the sense of Deligne (\cite[2.1.1]{De79}): \\[.1cm]
In fact, 
\cite{De79}  Proposition 1.1.14 derives the
axioms (2.1.1.1) and (2.1.1.2) from the fact that 
$\varphi_0$ defines a complex structure compatible with the
polarization. The axiom (2.1.1.3), i.e. the non-existence of a $\QQ$-factor
in $\Hg^{{\rm ad}}$ onto which $h$ projects trivially,
follows from $\Hg$ being the smallest $\QQ$-subgroup
of $\Sp(H^1(F,\QQ),Q)$ containing $\varphi_0$. By the condition
$(H_1)$ the bounded symmetric domains $X_{\Hg}$ and ${\mathcal X}$ 
coincide.
\par
\begin{Thm} \label{VHSShimura}
The VHS over a Shimura curve decomposes as follows:
$$ R^1 f_* \CC =: \VV_\CC = (\LL \otimes \TT) \oplus \UU$$
Here $\TT$ and $\UU$ are unitary local systems and $\LL$
is maximally Higgs, i.e.\ the corresponding Higgs field 
is an isomorphism. In particular, the VHS satisfies the 
Arakelov equality
\newline
If $C$ is non-compact then $\UU$ splits off over $\QQ$ and
over an unramified covering of $C$ the unitary local
systems become trivial.
\newline
Conversely if a family $f: \ol{A} \to \ol{C}$ of abelian
varieties satisfies the Arakelov equality, then its
VHS decomposes as above and $C$ is a Shimura curve.
\end{Thm}
\par
\begin{proof}
Since $\Hg$ is reductive we may split the representation
$\Hg \to \Sp(H^1(F,\QQ),Q)$ into a direct sum of irreducible
representations. We may split off unitary representations
over $\QQ$ (\cite{Ko87} Proposition~4.11). 
\par 
{\bf Claim:} For each of the
remaining irreducible representations there is an isogeny
$i: \SL_2(\RR) \times K \to \Hg_\RR$, where $K$ is a compact
group, such that the composition
of $i$ with $\Hg_\RR \to \Sp(H^1(F,\QQ),Q)_\RR$ is the tensor
product of a representation of $\SL_2(\RR)$ of weight one 
by a representation of $K$. 
\par
Assuming the Claim, consider a maximal compact
subgroup $K_1$ of $\SL_2(\RR) \times K$ that maps to the centralizer
of $\varphi_0$ under  
$$
\SL_2(\RR) \times K \longrightarrow \Sp(H^1(F,\QQ),Q)_\RR.
$$ 
The double quotient 
$$
X_{\Hg}' = \Gamma' \setminus\!(\SL_2(\RR) \times K)/K_1
$$ 
is an unramified cover of $X_\Hg$. Since a 
maximal Higgs field is characterized by the Arakelov
equality we may as well prove that the pullback
variation of Hodge structures has a maximal Higgs field. 
Since the fundamental
group of $X'_\Hg$ acts via $\Gamma'$ it follows immediately
from the claim and Lemma 2.1 in \cite{ViZu04} that the
VHS only consists of unitary and maximal Higgs subsystems.\par
\medskip\noindent
\par
{\it Proof of the Claim.}
We first analyze $\Hg^{{\rm ad}}$ and the $\Hg^{{\rm der}}(\RR)$-conjugacy
class $X_{\Hg}^{\rm ad}$ of maps ${\rm Res}_{\CC/\RR}\GG_m \to
\Hg^{{\rm ad}}$ containing $(\Hg \to \Hg^{\rm ad}) \circ \varphi_0$. 
Note that $X_{\Hg}^{\rm ad}$ is a connected component of $X_{\rm Hg}$. 
\newline
Each $\QQ$-factor of $\Hg^{{\rm ad}}$ onto which $h$ projects non-trivially
contributes to the dimension of $X_{\Hg}$. Since we
deal with Shimura curves $\Hg^{{\rm ad}}$ is $\QQ$-simple
by (2.1.1.3). Let 
$$\Hg^{{\rm ad}}_\RR = \prod_{i \in I} G_i$$
be its decomposition into simple factors. Then 
$X_{\Hg}^{\rm ad} = \prod X_i$ for $X_i$ a $G_i(\RR)$-conjugacy
class of maps ${\rm Res}_{\CC/\RR}\GG_m \to G_i$.
For the same reason, only one of the simple factors, say $G_1$, 
of $\Hg^{{\rm ad}}_\RR = \prod_{i \in I} G_i$
is non-compact. The possible complexifications $(G_1)_\CC$ are classified
by Dynkin diagrams. The property `Shimura curve', i.e.\ dimension one, implies 
that  $G_1 \cong {\rm PSl}(2,\RR)$.
\par
Now we determine the possible representations. The universal
cover $\widetilde{G_1} \to G_1$ factors though 
$$
H:={\rm Ker}(\Hg \to \prod_{i \in I \setminus \{1\}} G_i)^0.
$$ 
We apply \cite{De79} Section 1.3 to 
$$ (G_1,X_{\Hg}) \leftarrow (H,X_{\Hg}) \to (\Sp(H^1(F,\QQ)),\HH_g), $$
where $\HH_g$ is the Siegel half space and $g = \dim(F)$.
Since a finite-dimensional representation of $\widetilde{G_1}$ factors
through $\SL(2,\RR)$, we conclude that $G \cong \SL(2,\RR)$.
Moreover, such a representation corresponds to a fundamental
weight, hence of weight one. Now we let $\widetilde{K}$ be
the universal cover of $\prod_{i\in I\setminus \{1\}} G_i$. 
Since $\Hg \to \Hg^{\rm ad}$ is an isogeny, there is a lift 
of the universal cover $\widetilde{K} \to \Hg$. This lift
factors through a quotient $K$ of $\widetilde{K}$ such that the 
natural map $\SL(2,\RR) \times K \to \Hg_\RR$ is an isogeny. 
\par
Since we assumed the representation $\Hg_\RR \to \Sp(H^1(F,\QQ),Q)_\RR$
to be irreducible, also 
$$
\rho: \SL(2,\RR) \times K \to \Hg_\RR \to \Sp(H^1(F,\QQ),Q)_\RR
$$ 
is irreducible. Let $W \subset H^1(F,\RR)$ be an irreducible 
(necessarily weight one) representation of $\SL(2,\RR) \times \{{\rm id}\}$.
Since $K$ is reductive, hence its representations are semisimple, 
$\rho$ is the tensor product of $W$ and the representation
${\rm Hom}_{\SL(2,\RR) \times \{{\rm id}\}}(W,H^1(F,\RR))$ of $K$. 
This completes the proof of the claim.
\par
Since we have established the Arakelov equality, we may now refer
to \cite{ViZu04} to see that the VHS can be grouped together
as claimed. The statements in case $C$ non-compact are in \cite{ViZu04},
Section 4. 
\end{proof}
\par
The equivalence of the above definition of a Shimura curve and the
one stated at the beginning of the introduction is due to Moonen
(\cite{BMo98} Theorem~4.3).
\par
\begin{Rem} \label{re:trivialunitary}
If the base curve $C$ is not compact, then we may assume
that after passing to a finite unramified covering, the unitary
system $\UU$ is in fact trivial. This is shown in \cite{ViZu04}
Theorem~0.2.
\end{Rem}
\par
\section{Teichm\"uller curves} \label{Teich}
\par
\subsection{Moduli space of abelian differentials, $\SL_2(\RR)$-action}
\label{sl2rraction}
Let $\Omega M_g \to M_g$ be the moduli space of abelian differentials,
i.e.\ the vector bundle whose points are {\em flat surfaces}, i.e.\ 
pairs $(X, \omega)$ of a Riemann surface $X$ and a holomorphic $1$-form on $X$. 
The complement of the zero section, denoted by $\Omega M_g^*$, 
is stratified according to the number and multiplicities of
the zeros of $\omega$. This datum $(k_1,\ldots,k_r)$ with
$\sum k_i =2g-2$ is called the {\em signature} of $\omega$ or
of the stratum.
\par
There is a natural action of $\SL_2(\RR)$ on $\Omega M_g^*$ preserving
the stratification, defined as follows. Given a complex point
$(X,\omega)$ and $A \in \SL_2(\RR)$ we post-compose the
local charts of $X$ minus the set of zeros of $\omega$ 
given by integrating $\omega$ with the linear map $A$.
The composition maps give new complex charts of a Riemann surface
with punctures, that can be removed by Riemann's extension theorem.
We thus obtain a new compact Riemann surface $Y$.
The differentials $dz$ on the image of the charts inside $\RR^2 \cong
\CC$ glue to a well-defined holomorphic one-form $\eta$ on $Y$.
Consequently, we let $A\cdot(X,\omega)=(Y,\eta)$.
\par

\subsection{Teichm\"uller curves generated by abelian differentials}
\label{defTeich}
\begin{Defi}
A {\em Teichm\"uller curve} $C \subset M_g$
is an algebraic curve in the moduli 
space of curves, which is the image of the $\SL_2(\RR)$-orbit of
a flat surface $(X, \omega)$ under the projection $\Omega M_g \to
M_g$. In this case, the flat surface $(X, \omega)$ is called 
a {\em Veech surface}.
\end{Defi}
\par
In this situation we say that $C$ is {\em generated by} the 
Veech surface $(X, \omega)$. Teichm\"uller 
curves are algebraic curves in $M_g$ such that the universal
covering map $\HH \to T_g$ of the map $C \to M_g$ is a totally
geodesic embedding for the Teichm\"uller metric. In general, 
one can also define the $\SL_2(\RR)$-action on the space
of quadratic differentials over $M_g$ and include algebraic
curves generated by a pair $(X, q)$ with $q$ a quadratic
differential on $X$ in the definition of a Teichm\"uller curve, 
but we will not do this since we do not deal with these curves
in the sequel. 
\par
The normalization $C \to M_g$ of Teichm\"uller curve will still be
called Teichm\"uller curve. This is motivated by the following
 characterization.
\par
Let $\Aff^+(X,\omega)$ denote the group of orientation
preserving diffeomorphisms of $X$, that are affine with respect to
the charts determined by $\omega$. There is a well-defined
map $D: \Aff^+(X,\omega) \to \SL_2(\RR)$, that associates
with a diffeomorphism its matrix part. The image is called 
the {\em affine group} (or {\em Veech group}) of $(X,\omega$), 
usually denoted by $\Gamma$. Then the pair $(X, \omega)$ generates a
Teichm\"uller curve if and only if $\Gamma$ is a lattice
in $\SL_2(\RR)$. In this case $C = \HH/\Gamma$.

\subsection{Models of the universal family over a Teichm\"uller curve } 
\label{modelsUniv}

Parallel to the case of Shimura curves we consider the curves
$\HH/\Gamma_1$ for all $\Gamma_1 \subset \Gamma$ of finite index 
at the same time. For sufficiently small $\Gamma_1$ there
is a universal family of curves over $C_1:=\HH/\Gamma_1$ and the
monodromies around the cusps are unipotent, see the discussion
in Section~1 of \cite{Mo04a}. 
\par
Given a curve $C$ we denote by $\ol{C}$ its closure. For Teichm\"uller
curves $C$ and their unramified coverings defined above will drop the
subscripts from now on and -abusing the notations $C$, $\ol{C}$ etc.\ - 
assume that $\ol{f}: \ol{\cXX} \to \ol{C}$ is a stable family of curves. 
It is well-known that $f$ admits a model with smooth total space, 
if we relax the condition on the singular fibers from stable
to semi-stable. We will denote the unique relatively minimal model 
with semistable fibers and with smooth total space $\wtx$ by $\tilde{f}: 
\wtx \to \ol{C}$.
\par
The advantage of the stable model $\ol{\cXX}$ is its direct relation to the
geometry of Veech surfaces (see Section~\ref{sec:Degen}), whereas
calculations of intersection numbers work without correction terms
on $\wtx$ only.

\subsection{The characterization of Teichm\"uller curves by their VHS}
We recall from \cite{Mo04a} a characterization which, together with
Theorem~\ref{VHSShimura} will strongly restrict the kind of 
Teichm\"uller curves that may give rise to a ST-curve.
\par
\begin{Thm} \label{VHSTeich}
If $f: \cXX \to C$ is the universal family over an unramified covering
of a Teichm\"uller curve, 
then its VHS looks as follows: 
$$ R^1 f_* \CC =: \VV_\CC = (\oplus_{i=1}^r \LL_i) \oplus \MM. $$
Here $\LL_i$ are rank-$2$ subsystems, maximal Higgs for $i=1$ 
and non-unitary but not maximal Higgs for $i\neq 1$. In fact
$\MM$ splits off over $\QQ$ and $r=[K:\QQ]$, where $K=\QQ({\rm tr}(\gamma), 
\gamma \in \Gamma)$.
\newline
Conversely a family of curves $f: \cXX \to C$, whose VHS contains
a rank $2$ maximal Higgs subsystem, is the universal family
over a (covering of a) Teichm\"uller curve.
\end{Thm}
\par
\begin{Cor} \label{STimpliesr1}
If $C \to M_g$ is both a  Teichm\"uller curve and a Shimura curve, 
then $r=1$.
\end{Cor}
\begin{proof}
If $r>1$ then by Theorem~\ref{VHSTeich} the VHS over the Teichm\"uller
curve $C$ contains a local subsystem $\LL_i$ which is non-unitary but 
not maximal Higgs. Such a local subsystem does not appear in the decomposition of
the VHS over a Shimura curve by Theorem~\ref{VHSShimura}.
\end{proof}

\subsection{Some euclidian geometry}
\label{euclidian}
Given a Veech surface $(X, \omega)$. The differential
$\omega$ defines locally a euclidian coordinate system. With
respect to this the slope of a geodesic is well-defined. 
By abuse of notation we
call this slope and all geodesics with this slope a {\em direction}.
Veech dichotomy (\cite{Ve89}) states that each direction that contains
a geodesic joining two zeros or one zero to itself (a {\em saddle
connection}) is {\em periodic} i.e.\ each geodesic in this direction
is periodic or a saddle connection.
\par
A geodesic on a Veech surface $X$ has a {\em length} with respect
to the metric defined by $|\omega|$. The closed geodesics of a periodic direction (say the
horizontal one) sweep out cylinders
$C_i$ and we denote their core curves by $\gamma_i$. The {\em width}
of $C_i$ is defined to be the length of $\gamma_i$ in the metric $|\omega|$.
\par
Consider the degenerate fiber obtained by applying $\diag(e^t,e^{-t})$
to $(X, \omega)$ for $t \to \infty$. Say this point corresponds
to the cusp $c \in \ol{C} \sms C$. By \cite{Ma75} the
stable model of the fiber $\cXX_c$ of $f$ is obtained by squeezing
the core curves of the $C_i$ to points. Topologically the 
irreducible components of $\cXX_c$ are obtained cutting along the $\gamma_i$.
\par
 
\subsection{Square-tiled coverings} \label{sqtcov}
Teichm\"uller curves with affine group $\Gamma$ commensurable to $\SL_2(\ZZ)$, 
i.e.\ with $r=1$ arise as unramified coverings of the
once-punctured torus by \cite{GuJu00} Theorem~5.9. They are
called {\em square-tiled coverings} or {\em origamis}.
We denote torus covering by $\pi: X \to E$. Such a
covering can be uniquely factored as $\pi = i \circ \pi_{\rm opt}$, 
where $i$ is an isogeny and $\pi_{\rm opt}$ does not factor
through non-trivial isogenies. Such a $\pi_{\rm opt}$ is 
called {\em optimal} (sometimes also 'minimal' or 'maximal').
$\pi_{\rm opt}$ is ramified over torsion points. Given 
$\pi_{\rm opt}$ we may choose $\pi$ such that $i = [n]$ is
multiplication by some integer $n$.
\par
For a square-tiled covering the preimages
of torsion points are {\em periodic}, i.e.\ they have a finite 
$\Aff^{+}(X,\omega)$-orbit (\cite{GuJu00}, \cite{Mo04b}). Hence they
also define sections $R_i$ of $f$, after we maybe passed to a 
smaller $\Gamma$. The degeneration description also implies
that the $R_i$ do not intersect the sections $S_i$.
\par
For square-tiled surfaces
we normalize the length of horizontal and vertical saddle connections
to be integer-valued  by demanding that the horizontal and vertical closed geodesics
in $E \cong \CC/\ZZ[i]$ have length one.
\par

\subsection{Degeneration} \label{sec:Degen}
The following proposition 
shows that many boundary
components of the moduli space of stable curves $\ol{M_g}$
are not hit by the the closure of a Teichm\"uller curve in $\ol{M_g}$. 
\par
\begin{Prop} \label{noJsmoothsing}
The universal family over a Teichm\"uller curve $f:\cXX \to C$ does not contain singular fibers
whose Jacobian is proper, i.e.\ an abelian variety without toric
part. 
\end{Prop}
\par
\begin{proof} We give two proofs: The VHS depends only
on the family of Jacobians. If the Jacobian is proper for 
some $s \in S := \ol{C} \sms C$ the Higgs field 
$\Theta: E^{1,0} \to E^{0,1} \otimes \Omega^1_{\ol{C}}(\log S)$ 
factors through 
$E^{0,1} \otimes \Omega^1_{\ol{C}}(\log (S-s))$. This contradicts
that $\LL_1$ is maximal Higgs.
\newline
Or by topology: The Jacobian of a stable curve is proper if and only 
all the nodes are separating. As described in Section~\ref{euclidian} 
the nodes of the stable curves in the family over a Teichm\"uller
curves are obtained by squeezing the core curves of all
cylinders in a fixed direction. But a core curve of a cylinder can 
never be separating.
\end{proof}
\par

\subsection{Sections of the family of curves over a Teichm\"uller curve } \label{sectionsoff}
The zeros $Z_i$ of $\omega$ are periodic points (see \cite{Mo04b})
on the Veech surface. Hence they define sections of $f$, if we pass
to a sufficiently small subgroup of $\Gamma$ (see \cite{Mo04b} Lemma~1.2).
Let $S_i$ denote these sections. The unramified preimages under $\pi$
of $0 \in E$ are also periodic points. We will denote such a section
by $R$. From the above description of the degeneration we deduce that
neither the $S_i$ nor a section $R$ passes through a singularity of the
singulars fibers of the stable model $\ol{f}: \ol{\cXX} \to \ol{C}$. 
Consequently, the $S_i$ even define 
sections of $\ol{f}$ and we may speak of the zero $Z_i$ on a
singular fiber $\cXX_c$ of $\ol{f}$, meaning the intersection
of $S_i$ with $\cXX_c$.
\par

\section{Cyclic coverings} \label{sec:cyccov}

Families of cyclic coverings are easily understandable families of
curves with interesting variation of Hodge structures. They will
turn out to produce a number of Shimura and ST-curves. As we
will see in the next section, these are all examples with possible
exceptions in $g=5$.
\par
Fix integers $a_i$ and $N$ with $a_i \not\equiv 0 \mod N$, $a_1+a_2+a_3 \not\equiv 0 \mod N$.
For a family of cyclic coverings 
$$ \cYY_{N,\underline{a}}: y^N = x^{a_1}(x-1)^{a_2}(x-t)^{a_3} $$
the VHS decomposes into eigenspaces $\LL_i$ where the automorphism
$$\varphi: (x,y) \mapsto (x,\zeta_N y)$$
acts by the $i$-th power of $\zeta_N$.
\par
Their dimension and signature can be determined combinatorially
from $N$ and the 'type' $(a_1,a_2,a_3)$. This is known
since the work of Chevalley-Weil and is recalled for $4$-point covers
e.g.\ in Lemma~3.1 of \cite{BM05}.  The family  $\cYY_{N,(a_1,a_2,a_3)}$
defines a Shimura curve, if and only if on all the decomposition pieces
the monodromy action is by a Fuchsian group or by a unitary
group. Such a curve is moreover a Teichm\"uller curve, if
it contains no singular fiber with proper Jacobian. 
\par
For the special case $a_i=1$ for all $i$ we obtain the following
list of $ST$-curves.
\par
\begin{Thm} 
The family of curves $\cYY_{N,(1,1,1)}$ defines a ST-curve 
if and only if $N \in \{2,4,6\}$.
\end{Thm}
\par
\begin{proof}
Using the theory of admissible coverings one quickly shows
that the degenerate fibers for $t=0$ consists of two curves
isomorphic to $y^N=x(x-1)$, intersecting each other at $\gcd(N-2,N)$
normal crossings. Since degenerate fibers of the universal family
over a Teichm\"uller curve do not have separating nodes, this
implies for $ST$-curves that $N$ is even. 
\par
Since $a_1=a_2=a_3=1$, \cite{BM05} Lemma~3.1 implies that
for $i \leq N/3$ and for $i \geq 2N/3$ the Hodge filtration
on the local system $\LL_i$ is trivial, i.e.\ these local
systems are unitary. Since the local system $\LL_{N/2}$ 
pulls back from the family $y^2=x(x-1)(x-t)$ of elliptic
curves, it is maximal Higgs (with respect to
the set $\{0,1,\infty\}$ in the language of \cite{BM05}). 
This already implies that
for $N=4$ and $N=6$ the family is a Shimura curve. 
Since the set of singular fibers of the family $\cYY_{N,(1,1,1)}$
is exactly $\{0,1,\infty\}$, this family is a Teichm\"uller curve.
\par
It remains to show that for even $N>6$ one of the $\LL_i$
is neither maximal Higgs nor unitary. For this purpose
we take $i = N/2+1$. There is a general criterion
for checking if the monodromy group of a hypergeometric differential
equation (not only of second order) is unitary. Hypergeometric
differential equations of order $n$ are determined by a standard
set of parameters $\alpha_1,\ldots,\alpha_n$ and $\beta_1,\ldots,\beta_n$, 
see e.g.\ the introduction of \cite{BH89}. These $\alpha_i$ 
and $\beta_i$ are easily calculated from the local parameters
(compare the proof of \cite{BM05} Lemma~3.1) of the differential equation. 
These local parameters are, for the chosen $i$, given by $\{0,2/N\}$
at $x=0$, $\{0,2/N\}$ at $x=1$ and $\{(N-6)/2N,(N-2)/2N\}$. We obtain, since $N >6$,
$$0 < \alpha_1 = \frac{N-6}{2N} < \alpha_2 = \frac{N-2}{2N} < \beta_1 = \frac{2N-4}{2N}
< \beta_2 = 1.$$
The criterion of \cite{BH89} Corollary~4.7 is that $\exp(2\pi I (\cdot))$ applied to the
sets $\{\alpha_1,\alpha_2\}$ and $\{\beta_1,\beta_2\}$ interlace on the unit
circle. The above ordering shows that this is not the case. 
\par
In order to show that $\LL_{N/2+1}$ is not maximal Higgs we use
the uniqueness of a maximal Higgs subvariation of Hodge structures
of a family of curves, see \cite{Mo04a} Remark~2.5. 
\par
\end{proof}
\par
For the last argument one can, alternatively, calculate directly
the Lyapunov exponents (see Section~\ref{sec:Lyap}) of $\LL_{N/2+1}$ using
\cite{BM05} Proposition~3.4.
\par
For $N=4$ we next show that $\cYY_{4,(1,1,1)}$ defines a Teichm\"uller
curve by exhibiting the corresponding square-tiled covering. 
The same thing is easily done for $N=6$. The fact
that, at least for $N \geq 10$ the families $\cYY_{N,(1,1,1)}$
do not define ST-curves follows of course also from the main theorem,
which to be proven in the next section.
\par
For $N=5$ and $N=7$ the families $\cYY_{N,(1,1,1)}$ give
Shimura curves in the closure of $M_g$. This was first observed
in \cite{dJN91}. One can also use the methods of \cite{BM05}
Section~3 to show that the $\LL_i$ are either unitary or
maximal Higgs. Again, for odd $N > 7$, there exists a
subsystem $\LL_i$ that is neither maximal Higgs nor unitary.
\par

\subsection{The ST-curve in genus $3$}

We quickly present a ST-curve in genus $3$, as we will see the only
one, as a square-tiled covering. This way we can tell in the
next section that a remaining configuration is indeed this ST-curve. 
More details on this curve are in  \cite{HeSc04}.
\par
\begin{figure}[h] 
\centering
\epsfig{figure=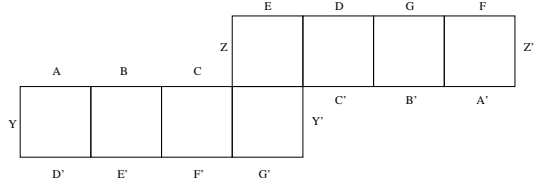, width=7cm}
\caption{The Teichm\"uller- and Shimura curve in genus $3$}
\label{figteichandshim}
\end{figure}

\begin{Prop}
The family of curves $ y^4 = x(x-1)(x-t)$ corresponds to the
square-tiled surface shown in Figure  \ref{figteichandshim}. 
The subgroup $\pi_1(E^*)$ defining
the covering is characteristic and has as quotient group the Quaternion
group $Q$ of order $8$. The Veech group of the square-tiled surface
is $\SL_2(\ZZ)$.
\end{Prop}
\par
\par
{\bf Proof:} Denote the covering map corresponding to the
above figure by $\pi: X \to E$. One checks that the
covering group is indeed $Q$.  In addition to that the elliptic involution 
lifts to an involution $\varphi$ of $X$. This can be checked
graphically since $180^\circ$ rotation of each square gives
a well-defined automorphism of $X$. Counting fixed points
one checks that $X/<\varphi>$ is of genus one. If $X$ was
hyperelliptic, the involution would have to be the lift of
the elliptic involution of $E$. Hence the generic fiber over 
the above Teichm\"uller curve is not hyperelliptic.
By \cite{KuKo79} there are only two families of curves of
genus three with an automorphism group of order $16$. The
one with hyperelliptic generic fiber can be ruled out and it
remains the above family.
\par
To determine the Veech group, remark that characteristic subgroups of $\pi_1(E^*)$ always have
Veech group  $\SL_2(\ZZ)$ (see e.g.\ \cite{He06}). 
\hfill $\Box$
\par

\section{Both Shimura- and Teichm\"uller curve}
\par
The proof of the following classification result will occupy the whole section:
\par
\begin{Thm} \label{Main}
For $g=2$ and $g \geq 6$ there is
no Shimura curve that is also a Teichm\"uller curve.
\newline
In both $M_3$ and in $M_4$ there is only one 
curve which is both a Shimura- and a {Teich}\-m\"uller curve. Its 
universal family is given by
$$ y^4 = x(x-1)(x-t) $$
in $M_3$ and respectively in $M_4$ by
$$ y^6 = x(x-1)(x-t). $$
In both cases  $t \in \PP^1 \sms \{0,1,\infty\}$.
\end{Thm}
\par
We recall that we abbreviate Shimura curves that are also 
Teichm\"uller curves as {\em ST-curves}. We now sketch the main steps
of the proof of this theorem. Existence statements were shown in 
the preceding Section~\ref{sec:cyccov}, so it remains to show the non-existence statements.
\par
{\em Strategy of the proof of Theorem~\ref{Main}:} 
From Corollary~\ref{STimpliesr1} and Section~\ref{sqtcov} it follows
that a ST-curve is necessarily a square-tiled covering, i.e.\ 
that $r=1$ or equivalently $\Gamma \subset \SL_2(\ZZ)$. Moreover, by
Theorem~\ref{VHSShimura} the local system $\MM$ appearing
in Theorem~\ref{VHSTeich} has to be
unitary. Since the ST-curve is not compact and by Remark~\ref{re:trivialunitary}, 
there is an abelian subvariety $A$ of dimension $g-1$ such
that $A \times C \hookrightarrow \Jac \wtx/\ol{C}$. The abelian variety $A$ is called
{\em fixed part} of the fibration. Fibered
surfaces with a fixed part of dimension $g-1$ have been called {\em maximally
irregular} in \cite{Mo02}. 
\par
In Section~\ref{degST} we show, for any $g \geq 3$, that all 
the singular fibers of the stable model $\ol{\cXX}$ are of a very particular type.
Using the dictionary between the cylinder decomposition of a Veech surface 
generating a ST-curve and the fibers over the cusps of the ST-curve, we
prove in Section~\ref{genus3} that for $g=3$ the ST-curve generated by the
flat surface in Figure~\ref{figteichandshim} is the only ST-curve. The proof has the advantage
to use flat geometry only, but Section~\ref{genus3} is strictly speaking
not necessary for a complete proof.
\par
In Section~\ref{Noether} we exploit the Noether formula (or equivalently
Riemann-Roch for $f: \cXX \to C$). This formula involves the self-intersection
of the relative dualizing sheaf of $f$. We are able to calculate this 
self-intersection (Lemma~\ref{omegarelrep} and Proposition~\ref{prop:less6}) 
exploiting two facts. 
By Section~\ref{degST} the family $f$
has no 'bad' singular fibers and moreover we know the self-intersections of 
sections of $f$ (Lemma~\ref{slefSi}) since $C$ arises as square-tiled surface. 
The positivity of the relative dualizing sheaf now implies the non-existence
for $g \geq 6$. 
\par
To obtain even more precise results for $g=4$ and $g=5$ we have to relate
the contributions of the singular fibers (the $\Delta \chi_{\rm top}(F)$ in
\eqref{eq:noether}) --  that we previously just estimated -- to the degree 
$d_\opt$ of the covering $\pi_\opt$ as defined in Section~\ref{sqtcov}.
\par
To prove the classification for $g=4$ we use that by \cite{Mo02} the
group $\SL(X,\omega)$ is contained the congruence subgroup $\Gamma(d)$
and that we may work over $\HH/\Gamma(d)$ all along. The geometric 
reformulation of this fact is given in Corollary~\ref{eqlenalldir}: We
may assume that $X$ is a branched cover of a square such that all cylinders
in all directions have the same widths. 
\par
On the other hand, from the Noether formula we obtain a quite restrictive list of
the singularities of a Veech surface generating a ST-curve for $g=4$ or $g=5$.
The list of singularities of course directly translates into the branching
behavior of $X \to E$. The remainder is a combinatorial discussion
whether such a covering with the cylinder condition from Corollary~\ref{eqlenalldir}
exists. \hfill $\Box$
\par

\subsection{Degeneration of ST-curves}\label{degST}
\par
We will only suppose $g>1$ here.
\par
\begin{Lemma} \label{geomgen}
The singular fibers of the smooth model of the universal family over a ST-curve $\tilde{f}$ 
have geometric genus $g-1$, 
hence $n$ components and $n$ nodes for some $n$. If we consider
the stable model $\ol{f}$ of the universal family, $n$ is at most $g-1$.
\end{Lemma}
\par
\begin{proof} The Jacobian of the singular fiber is a semiabelian
variety with non-trivial toric part by Proposition~\ref{noJsmoothsing}.
The abelian part has dimension equal to the geometric genus of the 
singular fiber.
Since the fixed part $A$ of dimension $g-1$ injects into the
Jacobian of each fiber, the geometric genus of the singular fiber
equals $g-1$. Hence the dimension of the toric part is one
and this dimension equals the rank of the fundamental group
of the dual intersection graph of the singular fiber. A graph
with fundamental group of rank one and without separating
nodes (compare the topological argument in the proof of 
Proposition~\ref{noJsmoothsing}) is a ring.
\end{proof}
\par
We now look at the cylinder decomposition of a fixed Veech surface
$(X,\omega)$ generating a ST-curve in a fixed 
periodic direction, say the horizontal one. 
As in Section~\ref{euclidian} let $C_i$ (for $i=1,\ldots,m$) 
denote the cylinders in this direction. 
\par
\begin{Lemma} \label{samelength}
For an appropriate numbering of the cylinders, the saddle connections bounding $C_i$ from the
above will bound $C_{i+1}$ from below (subscripts taken mod $n$).
Consequently, all cylinders of $(X,\omega)$ in a fixed direction 
have the same width.
\end{Lemma}
\par
\begin{proof}
By the previous lemma the number of components of the singular
fiber equals the number of nodes, which equals (again by
the description in \cite{Ma75}) 
the number of cylinders. By Proposition~\ref{noJsmoothsing}
this is only possible if the first
assertion is true and the second assertion is a consequence of the first.
\end{proof}
\par
In the sequel we use that we can retrieve singularities of $\omega$
also on the singular fibers, see Section~\ref{sectionsoff}.
\par
\begin{Lemma} \label{components}
If a component of a singular fiber of $\ol{f}$ contains zeros $Z_i$ 
of order $k_i$ for
$i=1,\ldots,\tilde{s}$, the geometric genus is $\sum_{i=1}^{\tilde{s}} k_i/2$.
In particular the geometric genus of a component of the
stable model $\ol{f}$ is always positive and
$\sum_{i=1}^{\tilde{s}} k_i$ is always even. 
\end{Lemma}
\par
\begin{proof} 
Gau\ss -Bonnet for a surface with $2$ boundary
components.
\end{proof}
\par
\begin{Lemma} \label{2cyl}
Each direction on a Veech surface generating a ST-curve contains at least two cylinders.
\end{Lemma}
\par
\begin{proof}
If not, we may find a transverse direction -- without loss of generality
the vertical one --  with a cylinder of width one. 
By Lemma~\ref{samelength} all cylinders in this
direction have this property. This means that the
square-tiled covering is abelian and $g(X)=1$.  
\end{proof}
\par
\begin{Lemma} \label{nosc}
A Veech surface generating a ST-curve has no saddle connection 
from a simple zero $Z$ to itself.
\end{Lemma}
\par
\begin{proof} 
Suppose that such a saddle connection $S$ connection exists, 
say on the top boundary of a horizontal cylinder $C_1$ and on the
lower boundary of $C_2$. Consider the
saddle connection $S'$ on the left of $S$, viewed from $C_1$. By
Lemma~\ref{samelength} it
must also be a saddle connection on $C_2$, whose right end is $Z$. 
But there are only two: $S$ and the saddle connection $S''$ 
on the left of $S$, viewed from $C_2$. If we had $S' = S''$ then
$Z$ would not be a zero at all. 
\end{proof}
\par
\par
\begin{Rem} {\rm
There are no ST-curves in genus two. This will
follow from Corollary~\ref{NoeCor}, but there are simpler proofs. E.g.\ 
\cite{Xi85} Theorem~3.16 implies that a genus two fibration
with one-dimensional fixed part has always singular fibers
with smooth Jacobian, a contradiction to Proposition~\ref{noJsmoothsing}.
}\end{Rem}
\par 
\subsection{Genus three} \label{genus3} \
\par
We suppose that $g=3$ and show that the square-tiled covering of
Figure~\ref{figteichandshim} generates the only ST-curve in $M_3$. We give here
a simple geometric proof not relying on the intersection
numbers on $\wtx$ as in Section~\ref{Noether}.
\par
Let $(X,\omega)$ be a Veech
surface generating a ST-curve in $M_3$. By Lemma~\ref{2cyl} and 
Lemma~\ref{geomgen} each periodic direction on $(X,\omega)$ 
has precisely two cylinders, and, by Lemma~\ref{geomgen} 
and Section~\ref{sec:Degen}, the corresponding singular fiber of the
stable model has precisely two components.
\par
\begin{Lemma} \label{g3allsimple}
All zeros of $\omega$ are simple.
\end{Lemma}
\par
\begin{proof} Suppose that there is a zero $Z$ of order three or four
and choose any periodic direction. The corresponding singular fiber of the
stable model has a component with  no zeros
or with only a simple zero. This contradicts Lemma~\ref{components}.
\par
Suppose there is a zero $Z$ of order two and choose a saddle connection
that connects $Z$ to some other zero $Z'$. Again the corresponding singular 
fiber of the stable model has a component containing $Z$ and $Z'$, while the
other component has at most one
simple zero. This contradicts again Lemma~\ref{components}.
\end{proof}
\par
\begin{Lemma} \label{sameall}
In each parabolic direction all saddle connections have the same length
and the cylinders have the same height.
\end{Lemma}
\par
\begin{proof}
Suppose we fix a parabolic direction where this is
not the case. Let $A$ be one of the longest saddle connections in this
direction. In particular $\ell(A) > w/4$, where $\ell$ denotes
the length of the saddle connection, normalized as in Section~\ref{Teich}, 
and $w$ the width
of any of the cylinders in this direction. By 
Lemma~\ref{g3allsimple} the square-tiled covering looks roughly 
as follows (see Figure~\ref{genus3rays}):
\par
\begin{figure}[h] 
\centering
\epsfig{figure=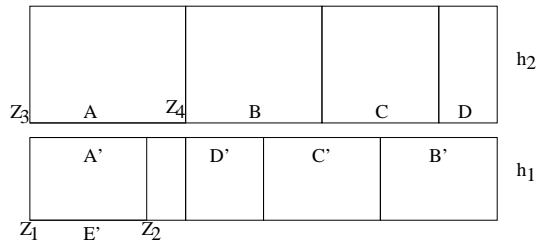, width=7cm}
\caption{A potential ST-curve in genus $3$}
\label{genus3rays}
\end{figure}
\par
The gluing $A$ to $A'$ etc.\ is imposed by the fact that we
have simple zeros and by Lemma~\ref{nosc}. The top line consists 
of four saddle 
connections $E$, $F$, $G$ and $H$, that are glued to
$E'$, $H'$, $G'$, $F'$ in this cyclic order on the bottom line.
We have only partially drawn the square tiling, the lengths of 
$A$, $B$, $E$, $F$, etc., $h_1$ and $h_2$ are integers. 
Turning the figure upside down we may suppose $h_1 \leq h_2$.
\par
Consider lines emanating from $Z_1$ passing through $A$. 
By Lemma \ref{nosc} the two occurrences of $Z_1$ on the
top line have to be in the segment not hit by such lines.
This segment is of length $w - \frac{h_1+h_2}{h_1} \ell(A)$.
It contains say $E$ and $F$. Hence 
$\ell(E) + \ell(F)$ equals at most
$w - \frac{h_1+h_2}{h_1} \ell(A)$. We may suppose
that 
$$\ell(E) \leq  \frac{1}{2}(w - \frac{h_1+h_2}{h_1} \ell(A)) < \ell(A),$$
where the strict inequality comes from the assumption 
that we want to contradict.
We suppose moreover that the 'short' saddle connection $E'$ is the one 
drawn in the Figure \ref{genus3rays}. 
For this we maybe have to choose the other
representative of $Z_1$ on the bottom line and/or look 'from
the back of the paper'.
\par
Now consider the lines emanating from $Z_2$ passing through $A$.
The segment on the top line hit by these lines does not contain
$Z_2$. The intersection $I$ with the corresponding 'hit segment'
from $Z_1$ does not contain any zero. One calculates
$$ \ell(I) = \frac{h_1+h_2}{h_1} \ell(A) - \frac{h_2}{h_1} \ell(E) 
 = \ell(A) + \frac{h_2}{h_1}( \ell(A) - \ell(E)) > \ell(A).$$
This contradicts the maximality of $A$.  
\end{proof}
\par
\begin{proof}[Proof of Theorem~\ref{Main} for $g=3$:] By the 
Lemmas~\ref{samelength}
and \ref{sameall} we know that the square-tiled covering is 
of degree $8$, maybe post-composed by unnecessary isogenies.
Since all zeros are simple, since there is no saddle connection
joining a zero to itself and since the, say, horizontal direction
has two cylinders of same length, the horizontal direction
consists of two cylinders of width $4$. There is only one
way to glue these cylinders to obtain simple singularities, 
the way shown in  Figure~\ref{figteichandshim}.
\end{proof}
\par
\subsection{Consequences of the Noether-formula on $\wtx$} \label{Noether}\ 
\par
We first show that an unramified covering of a ST-curve is modular
curve, i.e.\ the quotient of the upper half plane by a congruence subgroup. 
\par
\begin{Lemma} A ST-curve admits an unramified
covering $C=X(d)$, where $X(d)=\HH/\Gamma(d)$. Moreover
we can arrange that $d = \deg(\pi)$.
\end{Lemma}
\par
\begin{proof} By Theorem~1.6 of \cite{Mo02} a universal family of principally
polarized abelian varieties of dimension $g$ having as
fixed part some abelian variety 
$A$ of dimension $g-1$ exists over some $X(d_1)$. 
Hence the moduli map $C \to A_g$ factors through $X(d_1)$.
\par
What remains to show
is that in case the family comes from a ST-curve (i.e.\ consists
exclusively of Jacobians) the family of curves exists over
some $X(d)$, i.e.\ $C \to M_g$ factors through $X(d)$. 
\par
By loc.~cit.\  Proposition~1.7 the generic fiber of a ST-curve with
$g \geq 3$ cannot be hyperelliptic. We claim that $C \to A_g$
factors through $M_g$, if we take $\ol{C} \to \ol{X(d_1)}$ 
ramified of  order two at all the cusps. Indeed, the failure of 
infinitesimal Torelli at the hyperelliptic 
locus in $g \geq 3$ might force us to pass to a ramified cover, 
but a double cover is always sufficient. Moreover the cover does not
need to be ramified inside $X(d_1)$ since otherwise the local system
$\LL_1$ of the Teichm\"uller curve would no longer be
maximal Higgs but acquire a zero at the ramification point. 
For large enough  $d/d_1$ the covering
$\ol{X(d)} \to  \ol{X(d_1)}$ is ramified at all cusps of $\ol{X(d_1)}$
of even order. Hence the ST-curve exists over $X(d)$ for some $d$.  
\par
The second statement can be arranged by post-composing $\pi$ with
an isogeny.
\end{proof}
\par
In the sequel we will sometimes enlarge $d$ further if necessary: We
suppose that $d$ is large enough such that the base $X(d)$ has positive genus.
Since the fibers have genus greater than two, $\cXX$ has a unique minimal
model and this hypothesis on $X(d)$ ensures that the smooth relatively minimal model 
for $f:X \to C$ is in fact minimal.
\par
On the smooth minimal model  $\tilde{f}:\wtx \to \ol{C}$ we have the 
Noether equality 
$$ 12 \chi(\cOO_{\wtx}) - c_2(\omega_{\wtx}) = c_1(\omega_{\wtx})^2.$$
For a fibered surface with fiber genus $g$ and base genus $b$ we have
$$ \begin{array}{lcl}
\chi(\cOO_{\wtx}) & =& \deg \tilde{f}_* \omega_{\wtx/\ol{C}} + (g-1)(b-1)\\
 c_2(\omega_{\wtx}) &=& \sum_{F\;{\rm sing.}} 
\Delta \chi_{\rm top}(F)+ 4(g-1)(b-1)\\
c_1(\omega_{\wtx})^2 &=& \omega^2_{\wtx/\ol{C}} + 8 (g-1)(b-1).\\
\end{array}$$
We will hence use the equality in the form 
\begin{equation} \label{eq:noether}
12 \deg \tilde{f}_* \omega_{\wtx/\ol{C}} - \sum_{F\;{\rm sing.}}
\Delta \chi_{\rm top}(F)
 = \omega^2_{\wtx/\ol{C}}. 
\end{equation}
Here the sum runs over all singular fibers $F$ of $\tilde f$ and 
$\Delta \chi_{\rm top}$ denotes the difference of the topological
Euler characteristic of $F$ and $2-2g$. 
\par
It is well-known (e.g.\ \cite{Sh71}) that
$$ g(X(d)) = (d-6)\Delta_d + 1,\, \text{where}\, \Delta_d = \left\{
\begin{array}{ll} \frac{d^2}{24}
\prod_{p|d} (1-\frac{1}{p^2}) & d \geq 3 \\
\frac{1}{4} & d = 2 \end{array} \right.$$
and that the number of cusps is
$$ |S| := |\ol{X(d)} \sms X(d)| = 12 \Delta_d, \quad \text{hence} \quad
\deg \Omega^1_{\ol{C}}(\log S) = 2 d \Delta_d.$$
For a Shimura curve with fixed part of dimension $g-1$ the Arakelov equality 
states that 
$$\deg \tilde{f}_* \omega_{\wtx/\ol{C}} = \frac{1}{2}\deg \Omega^1_{\ol{C}}(\log S)
=  d \Delta_d.$$
\par
We start calculating $c_1(\omega_{\wtx})^2$ and $c_2(\omega_{\wtx})$
for a square-tiled covering. Note that the map $\varphi: \wtx \to \ol{\cXX}$ 
from the smooth semistable model to the stable model 
only contracts $(-2)$-curves. Hence $\omega_{\wtx/\ol{C}} = 
\varphi^* \omega_{\ol{X}/\ol{C}}$.
Let $h: \cEE \to X(d)=C$ be the universal family of elliptic curves 
with full level-$d$-structure. We extend this family to a smooth
minimal model $\tilde{h}: \widetilde{\cEE} \to \ol{C}$.
\par
Recall from Section~\ref{sectionsoff} the 
the definition of the sections $S_i$ and $R$. We need to calculate
the self-intersection numbers of these sections. In order to do
so we have to perform some base changes and blowups, that will
not appear in the formulae at the end of the day.
\par
The sections $S_i$ and $R$ are defined over 
some unramified cover $\psi:C_{\rm sec} \to X(d)$, not necessarily
given by a congruence subgroup. Let 
$$\psi_X: \widetilde{\cXX} \times_{\ol{C}} \ol{C}_{\rm sec}
\to \widetilde{\cXX}$$ 
be the covering of the surfaces induced by the base change $\psi$. 
There is a blowup
$b: Y \to \widetilde{\cXX} \times_{\ol{C}} \ol{C}_{\rm sec}$ 
in the singularities of the 
singular fibers of $\widetilde{\cXX} \times_{\ol{C}} \ol{C}_{\rm sec} 
\to \ol{C}_{\rm sec}$
such that $Y$ is smooth and such that the square-tiled 
covering map $\pi$ extends to a map
$$\pi_{\rm sec}: Y \to \widetilde{\cEE_{\rm sec}}:= 
\widetilde{\cEE} \times_{\ol{C}} \ol{C}_{\rm sec}. $$
By definition of $Y$ there is a blowdown $b_2:Y \to \widetilde{\cXX_{\rm sec}}$
to the minimal semistable model $\widetilde{\cXX_{\rm sec}}$ over 
$\ol{C}_{\rm sec}$, through which $b$ factors.
\par
Suppose that the Teichm\"uller curve lies in the stratum with
signature $(k_1,\ldots,k_s)$, i.e.\ in each smooth fiber $F$ of $f$ 
the zeros of $\omega|_F$ have multiplicity $k_1,\ldots,k_s$. 
\par
\begin{Lemma} \label{slefSi}
The sections $S_i$ on $\widetilde{X_{\rm sec}}$ have self-intersection 
number 
$$S_i^2 = \frac{-\deg \tilde{f}_* \omega_{\widetilde{X_{\rm sec}}/\ol{C_{\rm sec}}}}
{k_i+1}.$$
\end{Lemma}
\par
\begin{proof} Kodaira has calculated (\cite{Kd63} (12.6) and (12.7)) 
that the zero-section $\sigma_0$ of $\cEE$ has
self-intersection 
$$\sigma_0^2 = -d\Delta_d = -\deg \tilde{f}_* \omega_{\widetilde{\cXX}/\ol{C}},$$ 
hence on the pullback surface $\widetilde{\cEE_{\rm sec}}
:= \widetilde{\cEE} \times_{\ol{C}} \ol{C_{\rm sec}}  
\to \ol{C_{\rm sec}}$ we have
$$(g_{\cEE}^* \sigma_0)^2 = -\deg \tilde{f}_* \omega_{\widetilde{\cXX_{\rm sec}}/
\ol{C_{\rm sec}}},$$
where $g_{\cEE}: \widetilde{\cEE_{\rm sec}} \to \widetilde{\cEE}$.
The blowup $b_2$ is an
isomorphism outside the singularities of the singular fibers of $f$
and $S_i$ does not hit them, see the description of the degeneration
in Section \ref{Teich}.
Hence $$S_i^2 = (b_2^* S_i)^2 = \frac{1}{k_i+1} (b_2^* S_i \cdot \pi^* g_\cEE^*
\sigma_0) = \frac{-\deg f_* \omega_{\widetilde{X_{\rm sec}}/\ol{C}_{\rm sec}}}
{(k_i+1)}.  $$  
\end{proof}
\par
The same argument applies to the section induced by an unramified 
torsion point 
$R$ and yields $R^2 = -d\Delta_d$. In the sequel we denote the families of
curves over $\ol{C}_{\rm sec}$ still by $\ol{f}$ or $\tilde{f}$, without
the subscript ${\rm sec}$.
\par
\par
\begin{Lemma} \label{omegarelrep}
On the smooth minimal model $\tilde{f}:
\wtx_{\rm sec} \to \ol{C}_{\rm sec}$ of a ST-curve coming
from a square-tiled covering we have 
$$ \omega_{\wtx_{\rm sec}/\ol{C}_{\rm sec}} = {\cOO}_{\wtx_{\rm sec}}
\left(\sum_{i=1}^s k_i S_i +  
\deg \tilde{f}_* \omega_{\widetilde{X_{\rm sec}}/\ol{C}_{\rm sec}} \cdot F\right), $$
where $F$ denotes a fiber.
\end{Lemma}
\par
\begin{proof}
The generating differential $\omega^0$ defines a subbundle $\cLL \subset
\tilde{f}_* \omega_{\wtx_{\rm sec}/\ol{C}_{\rm sec}}$. By definition of
the $k_i$ and $S_i$, the cokernel of
$\tilde{f}^* \cLL(\sum k_i S_i) \to \omega_{\wtx_{\rm sec}/\ol{C}_{\rm sec}}$ 
is supported in the singular fibers.
\par
We claim that the cokernel does not contain components of singular
fibers but only the complete fibers. Once the claim is established
the number of these fibers is $-R^2= \deg \tilde{f}_* \omega_{\widetilde{X_{\rm
sec}}/\ol{C}_{\rm sec}}$ by the adjunction formula.
\par
To prove the claim, suppose the cokernel contains a connected
divisor $B$ with multiplicity $\mu >0$, contained in a singular fiber 
but not equal to
a singular fiber. Lemma~\ref{components} also applies to $B$, 
since the dual graph of $B$ contains no loops and yields
$g(B) = \sum_{i\in I_B} k_i/2$, where the sum is over the sections $S_i$,
$i \in I_B$
that hit the component $B$. On the other hand $B^2 = -2$, since 
$B$ hits the rest of the fiber in two points by Lemma~\ref{geomgen}. 
Now the adjunction formula says
$$ B \cdot \omega_{\wtx_{\rm sec}/\ol{C}_{\rm sec}} = 2g(B)-2 - B^2 = 
2g(B) = \sum_{i \in I_B} k_i.$$
On the other hand, by our hypothesis on the cokernel
$$  B \cdot \omega_{\wtx_{\rm sec}/\ol{C}_{\rm sec}} = \mu B^2 +
\sum_{i \in I_B} k_i
$$
This is a contradiction.
\end{proof}
\par
Although we will calculate more precisely, let us note here
an immediate corollary
\par
\begin{Prop} \label{prop:less6}
ST-curves exist only for $g \leq 6$.
\end{Prop}
\par
\begin{proof} Since $g \geq 2$ and since we may suppose that
$g(C_{\rm sec}) \geq 2$ by enlarging $d$, the surface $\wtx$ is
of general type and hence $\omega_{\wtx_{\rm sec}/\ol{C}_{\rm sec}}$ is nef.
Since there is at least one singular fiber we have
$$ 
\begin{array}{lcl} 
12 d \Delta_d  &>&
(\omega_{\wtx_{\rm sec}/\ol{C}_{\rm sec}} \cdot d \Delta_d F) + 
(\omega_{\wtx_{\rm sec}/\ol{C}_{\rm sec}}
\cdot (\sum_{i=1}^s k_i S_i)) \\
& \geq & (\omega_{\wtx_{\rm sec}/\ol{C}_{\rm sec}}
\cdot d \Delta_d F) = (2g-2)d\Delta_d.
\end{array}$$
Solving for $g$ gives the desired inequality.
\end{proof}
\par
We now come to the calculation of $\Delta \chi_{\rm top}$. For this
purpose, it will be more convenient to work with family of curves 
over $C = X(d)$ instead of the covering over $C_{\rm sec} \to C$. 
Recall that from Section~\ref{sqtcov} that we may 
suppose $\pi = [n] \circ \pi_{\rm opt}$. With this normalization,
the image of $\pi_{\rm opt}$ is the universal family of elliptic 
curves over $X(d)$, not just any family isogenous to that one.
\par
\begin{Prop} \label{Neronprop}
For the smooth minimal model $\wtx$ of a ST-curve we
have for each of the $12\Delta_d$ singular fibers 
$$\Delta \chi_{\rm top} = d/d_{{\rm opt}}.$$
\end{Prop}
\par
\begin{proof} 
By Lemma \ref{geomgen} the singular fibers of $\tilde{f}$
consist of a ring of say $\tilde{m}$ smooth curves. In that
case $\Delta \chi_{\rm top} = \tilde{m}$. Fix some 
$c \in S = \ol{C} \sms X(d)$ and let $\ol{U}$ be a neighborhood
of the cusp, i.e.\ the spectrum of a discrete valuation ring
whose closed point is $c$. The N\'eron-Model of the
Jacobian $J_U$ of the family of curves $\tilde{f}|_U$ is 
$J_{\ol{U}} := {\rm Pic}^0_{\widetilde{\cXX}|_{\ol{U}}/\ol{U}}$.
The special fiber $J_c$ of $J_{\ol{U}}$ is an an extension of
a semiabelian variety by a finite group, the local component 
group $\Phi_{J_c}$. By \cite{BLR90} Remark~9.6.12 the cardinality of
$\Phi_{J_c}$ equals  $\tilde{m}$.
\par
We want to compare this to the N\'eron-Model $E_{\ol{U}}$ of the
Jacobian $E_U$ of the family of elliptic curves $h: \cEE \to \ol{C}$
restricted to $U$. The component group $\Phi_{E_c}$ of the
special fiber of $E_{\ol{U}}$ is known to have order $d$, since
we work over $\ol{C}=\ol{X(d)}$.
\par 
Recall that $A$ is the fixed part in the family of Jacobians over $C$.
We let $A_{\ol{U}}:= A|_U \times_U \ol{U}$. Then the 
exact sequence of abelian varieties
$$ 0 \to A|_U \to J_U \to E_U\to 0$$
induces by \cite{BLR90} Proposition~7.5.3 a) an exact sequence
$$ 0 \to A_{\ol{U}} \to J_{\ol{U}} \to E_{\ol{U}}$$
of N\'eron models. The right arrow will no longer be
surjective.
\par
Since $A$ is constant, we deduce from this sequence an injection 
$\Phi(\pi_{{\rm opt}}): \Phi_{J_c} \to \Phi_{E_c}$.
The cokernel of $\Phi(\pi_{{\rm opt}})$ has been determined
in \cite{CoSt01} Theorem~6.1: To check their hypothesis, note 
that $\pi_{{\rm opt}}$ is optimal by definition, the degeneration 
of $E_U$ is toric and Theorem~6.1 in loc.\ cit.\  
not only works over a p-adic field 
but over any field complete with respect to a non-archimedian
valuation. We introduce some notation to use this result 
in the form of loc.~cit.\ Corollary~6.6.
Let $X_{E_c}$ resp.\ $X_{J_c}$ be the character group of the
toric part of $E_c$ resp.\ of $J_c$. There is a natural
map $\pi_{{\rm opt}}^* X_{\cEE_c} \to X_{J_c}$ and we let
$\cLL$ be the saturation of the image. The first quantity
we need is $m_{\cLL} = [\alpha(X_{J_c}): \alpha(\cLL)]$.
The homomorphism $\alpha$ is induced by the monodromy 
pairing. We do not recall details, since in our special situation 
we need not any property of $\alpha$ other than being a homomorphism.
Second, let $\Theta$ be the polarization on $E_U$ obtained 
by pulling back the principal polarization on $J_u$ via $\Jac(\pi_{\rm opt})$.
We define $m_{E_U} = \sqrt{\deg(\Theta)}$. Now the formula 
of \cite{CoSt01} for the size of the cokernel is 
$$ |{\rm coker}(\Phi(\pi_{{\rm opt}}))| = 
m_{E_U} / m_{\mathcal L}$$
\par
In our special case, the intersection graphs $\Gamma_{\cEE_c}$ and
$\Gamma_{J_c}$ of the components of the singular fibers of both
$\cEE_c$ and $J_c$ are rings by Lemma~\ref{geomgen}. 
Hence both $H_1(\Gamma_{\cEE_c},\ZZ)\cong 
X_{\cEE_c}$ and  $H_1(\Gamma_{J_c},\ZZ) \cong X_{J_c}$ are cyclic. 
We conclude that $\cLL = X_{J_c}$ and thus $m_{\cLL} =1$.
\par
On the other hand, $m_{E_U} = d_{\rm opt}$, see e.g.\ \cite{BiLa03} 
Lemma~12.3.1. Assembling all the data, we obtain 
$$\Delta \chi_{\rm top}= \tilde{m} = d /|{\rm coker}(\Phi(\pi_{{\rm opt}}))|
= d/d_{\rm opt}, $$ 
as claimed.
\end{proof}
\par
When we insert the results of the previous lemmas into the
Noether formula we obtain:
\par
\begin{Cor} \label{NoeCor}
For a ST-curve we have 
$$ d_{{\rm opt}} = \frac{12}{\sum_{i=1}^s \frac{k_i^2}{k_i+1} + 16-4g}.$$
Hence such curves exist at most for $g=3$, signature $(1,1,1,1)$
and $d_{{\rm opt}} = 2 $ -- which have been studied in 
Section \ref{genus3} -- or for one of the
cases in the following table:
\newline
\centerline{
\begin{tabular}{cccc}
case & genus & signature & $d_{opt}$ \\
(1)    & 4     & (1,1,1,1,1,1) & 4 \\
(2)    & 4     & (2,2,2)       & 3 \\
(3)    & 5     & (1,1,1,1,1,1,2) & 36 \\
(4)    & 5     & (1,1,1,1,2,2) & 18 \\
(5)    & 5     & (1,1,2,2,2) & 12 \\
(6)    & 5     & (2,2,2,2) & 9 \\
(7)    & 5     & (1,1,1,1,1,3) & 16 \\
(8)    & 5     & (1,1,3,3) & 8 \\
(9)    & 5     & (1,1,1,1,4) & 10 \\
\end{tabular}}
\end{Cor}
\par
\begin{proof} From Lemma \ref{slefSi} and Lemma \ref{omegarelrep}
we deduce that on $\wtx_{\rm sec}$
$$  (\omega_{\wtx_{\rm sec}/\ol{C}_{\rm sec}})^2
= \left(2(2g-2) - \sum \frac{k_i^2}{k_i+1}\right) 
\deg \tilde{f}_* \omega_{\widetilde{X_{\rm sec}}/\ol{C}_{\rm sec}}. $$
Since $\deg \tilde{f}_* \omega_{\widetilde{X_{\rm sec}}/\ol{C}_{\rm sec}}
= \deg(\psi) \deg \tilde{f}_* \omega_{\widetilde{X}/\ol{C}}
= \deg(\psi) \cdot d \Delta_d$, we have on
$\wtx$
$$ (\omega_{\wtx/\ol{C}})^2 = \left(2(2g-2) - \sum \frac{k_i^2}{k_i+1}
\right) \cdot d \Delta_d. $$
Together with Proposition \ref{Neronprop} and the Noether formula this yields
the claimed expression for $d_{\rm opt}$.
\par
We now list all the signatures that provide non-negative integer solutions of the
equation for $d_{\rm opt}$, but which are missing in the
table: $(4,4)$ and $(1,2,5)$ with $g=5$ and $(3,7)$, $(1,1,8)$, $(2,8)$, 
$(1,9)$, $(5,5)$ and $(10)$
with $g=6$. We argue why they cannot occur for ST-curves.
\par
In the cases with two or less zeros, consider a direction 
containing a saddle connection joining the two zeros. By
Lemma~\ref{geomgen} the stable curve obtained from degenerating in this 
direction has only one component.
Hence there is only one cylinder in this direction by Lemma~\ref{components}.
This contradicts Lemma~\ref{2cyl}.
\par
In the remaining cases, take a direction containing a saddle 
connection joining the zero of order $5$ (resp.\ $8$) with the
zero of order $2$ (resp.\ a zero of order $1$). When degenerating
in this direction, the remaining zero lies by Lemma~\ref{components} 
on the same component as the other two zeros. We now conclude as
in the above cases.
\end{proof}
\par
Meanwhile we can significantly strengthen Lemma~\ref{samelength}.
\par
\begin{Cor} \label{eqlenalldir}
With the normalization that $\pi = [n] \circ \pi_{{\rm opt}}$ the
width $w_\theta$ of each of the cylinders is the same in each 
direction $\theta$. In particular, the number of $\pi_{\rm opt}$-preimages
of a cylinder in the direction $\theta$ in $E_n$ is independent 
of the cylinder and independent of $\theta$.
\end{Cor}
\par
\begin{proof}
Suppose that in a fixed direction $\theta$ the ST-curve has $k_\theta$
cylinders $C_i$ of height $h_i$. Consequently, $d= (\sum_{i=1}^{k_\theta} h_i) 
\cdot w_\theta$ where
$w_\theta$ is the width of the cylinder in the direction $\theta$.
Let $c = c(\theta)$ be the cusp corresponding to the degeneration
in the direction $\theta$. We know by Proposition~\ref{Neronprop} 
and the normalization of $\pi$ that the semistable model $\wtx$
of the fiber over $c$  is a ring of $\tilde{m}$
curves, independently of the direction. The stable fiber
over $c$ has $k_\theta$ components and we have to determine 
the type of the singularities. The fundamental group locally
acts by $\gamma = \left(\begin{array}{cc} 1 & d \\ 0 & 1\\ \end{array}\right)$ 
since we work over the base $\HH/\Gamma(d)$.
Thus the cylinder $C_i$ is twisted $dh_i/w_\theta$ times 
by $\gamma$. This leads to
a singularity of type $A_{dh_i/w_\theta}$ at the corresponding node of
the singular fiber. After resolving by a string of $dh_i/w_\theta-1$ projective
lines the semistable singular fiber is a ring of 
$$ \tilde{m} = \sum_{i=1}^{k_\theta} \left(h_i \frac{d}{w_\theta} -1\right) + {k_\theta} = 
\left(\frac{d}{w_\theta}\right)^2$$
curves.
\end{proof}
\par 
We consider now the flat surface $X$ as a branched 
covering of $E_n = \CC/(n\ZZ[i])$
with the degree and ramification type being one of the cases in
Corollary~\ref{NoeCor}. We denote by $B \subset E_n$ the
set of branch points, we let  $E_n^* := E_n \sms B$ and
apply  Corollary~\ref{eqlenalldir} to the decomposition
of $E_n^*$ into cylinders in various directions.
\par
\begin{Lemma} \label{onecylpreim}
In each of the possible cases of ST-curves listed in Corollary \ref{NoeCor} 
one of the following two possibilities holds.
\begin{itemize}
\item[i)] The preimage under $\pi_{{\rm opt}}$ of each cylinder in $E_n$ 
consists of only one cylinder in $X$, or 
\item[ii)] $B$ consists of one element only and each cylinder in $E_n$ 
has $k$ preimages under $\pi_{\rm opt}$ where $2 \leq k \leq g-1$
, or
\item[iii)] $g=5$ and $B$ is the set of $2$-torsion points of $E_n$.
Moreover, each cylinder in $E_n$ has $2$ preimages under $\pi_{\rm opt}$ 
\end{itemize}
\end{Lemma}
\par
\begin{proof} Suppose first that $g=4$. Suppose there exists 
some direction $\theta$ such that $E_n^*$
decomposes into more than one cylinder in the direction $\theta$.
If one cylinder in the direction $\theta$ in $E_n$ has more than 
one preimage under $\pi_{\rm opt}$, all cylinders in the direction $\theta$
have this property by Corollary~\ref{eqlenalldir}. Since the number of 
horizontal cylinders in $X$ is bounded by $g-1 =3$ by Lemma~\ref{geomgen} 
and Lemma~\ref{components}, we obtain a contradiction.
 Corollary~\ref{eqlenalldir} now implies the same behavior in
all the other directions
\par
The non-existence of a direction $\theta$ as above implies
that $B$ consists of one element only. The upper bound on $k$ in this
case is obtained from the same argument. The lower bound is
Lemma~\ref{2cyl}.
\par
In the case $g=5$ we can apply the same argument with one exception:
Each cylinder in $E_n$ has precisely two preimages in $C$ and
each direction of $E_n$ consists of at most two cylinders. The
last condition implies that the set of branch points of $\pi_{{\rm opt}}$
equals the set of $2$-torsion points of $E_n$. 
\end{proof}
\par
We let $R$ be the ramification points of $\pi_{\rm opt}$ and we
present the unramified covering
$$\pi_{{\rm opt}}: X^* := X \sms R \to E_n^*$$
by a homomorphism $\varphi_\opt: \pi_1(E_n^*, P) \to S_{d_{{\rm opt}}}$, where
$P$ is some base point. 
\par

\subsection{Case $d_{{\rm opt}}$ is even} 

The condition in Lemma~\ref{onecylpreim}~i) 
implies that elements of $\pi_1(E^*, P)$ representing core curves
of cylinders are mapped to $d_{{\rm opt}}$-cycles. 
The hypothesis $d_{{\rm opt}}$ is even translates into 
a $d_{{\rm opt}}$-cycle being odd and we thus can argue with parities. 
We illustrate the basic idea.
Suppose that $\alpha$ (resp.\ $\beta$) represents the core curve
of a horizontal (resp.\ vertical) cylinder. Then 
$\varphi_\opt(\alpha)$ and $\varphi_\opt(\beta)$ are odd. If none
of the branch points was mapped to an odd element in $S_{d_{\opt}}$,
then the core curve of a cylinder of slope $1$ would have the same
parity as $\varphi_\opt(\alpha)\varphi_\opt(\beta)$, that is even parity.
It can thus not correspond to a $d_{{\rm opt}}$-cycle and we obtain
a contradiction.
\par
We call branch points of $\pi_\opt$ whose monodromy is even (resp.\ odd) 
in $S_{d_\opt}$ {\em even} (resp.\ {\em odd}) {\em branch points}. 
They must be in a very special
position in order to avoid a contradiction similar to the above one.
\par
\begin{Lemma} \label{le:dopteven}
If  $d_{{\rm opt}}$ even and if case i) of Lemma~\ref{onecylpreim} holds, 
then it is possible to normalize $E_n$ by translating the origin such
that the set of odd branch points is precisely the set of 
$2$-torsion points of $E_n$.
\end{Lemma}
\par
\begin{proof} The set of odd branch points is non-empty by
the argument preceding the Lemma and by translation one of them is
the origin. By the genus condition of Lemma \ref{components}
each straight line of rational slope $E_n$ contains an even number
of odd branch points. By the list of Corollary~\ref{NoeCor} there are 
at most $6$ odd branch points. From these two statements we
deduce that there are precisely $4$ odd branch points,
given by $P_1 = (0,0)$,  $P_2= (0,a)$, $P_3 = (a,0)$
and $P_4 = (a,a)$ for some $a \in \{0,\ldots,n/2\}$.
\par
If $a \neq n/2$, the diagonal of slope $-1$ contains only
one odd branch point, namely $P_1$. This is a
contradiction and completes the claim.
\end{proof}
\par
\begin{Cor} \label{cor:nocase1}
There are no  ST-curves in the case (1). In the cases with $g=5$
and $d_\opt$ even (i.e.\ cases (3),(4),(5),(7),(8),(9)) there is only a finite list of candidates
of square-tiled coverings that could give rise to a ST-curve.
\end{Cor}
\par
\begin{proof}
In case (1), the case ii) of Lemma~\ref{onecylpreim} cannot appear, since in this
case the unique branch point had $12 > d_\opt = 4$ preimages. In all the
cases with $g=5$ the possibilities ii) and iii)  Lemma~\ref{onecylpreim} give
a finite number of cases only, since the degree and the location of
the branch point(s) are known.
\end{proof}
\par
We remark that the preceding argument gives another proof of
the classification in genus three.
\par

\subsection {Case $g=4$, signature $(2,2,2)$, and $ d_{{\rm opt}} = 3$}
We represent the covering by its monodromy $\varphi_\opt$
as in the previous section. 
$ d_{{\rm opt}} = 3$ implies that there are in fact
precisely three branch points. By applying a matrix in
$\SL_2(\RR)$ we may assume that the branch points are
located at $P_0=(0,0)$, $P_1=(a,0)$ and $P_2=(0,b)$.
\begin{figure}[h] 
\centering
\epsfig{figure=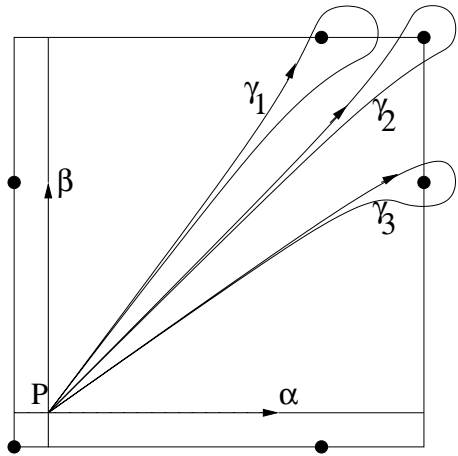, width=6cm}
\caption{Standard presentation of $\pi_1(E_n^*)$ punctured at three
points.}
\label{fig:3punct}
\end{figure}
\par
\begin{Lemma} \label{le:case2}
The only ST-curve in the case (2) is
given by $y^6=x(x-1)(x-t)$.
\end{Lemma}
\par
{\bf Proof:}
Consider the images of $\alpha$, $\beta$ and the $\gamma_i$
under the monodromy map $\varphi_{{\rm opt}}$.  
Lemma~\ref{onecylpreim} and the hypothesis on the signature 
imply respectively that they all map to $3$-cycles. 
Since these paths generate $\pi_1(E_n^*,P)$, the map $\varphi_\opt$
factors through $\ZZ/3 = \langle \sigma \rangle \subset S_3$.
Since $[\alpha,\beta]=\gamma_1\gamma_2\gamma_3$ we must have
$\varphi(\gamma_i) =\sigma$ for $i=1,2,3$, choosing $\sigma$
suitably. 
\par
By symmetry, we may suppose that $a \geq n/2$ and $b \geq n/2$.
Since $\alpha$ and  $\alpha\gamma_3$ represent core curves of cylinders, 
we deduce $\varphi(\alpha) =\sigma$. The same argument applied
to $\beta$ and $\beta \gamma_1^{-1}$ yields $\varphi(\beta) =\sigma^2$.
We thus have determined uniquely the branched cover up to
the location of $a$ and $b$. If $a=b=n/2$, e.g.\ a calculation
of the automorphism group shows that the branched
cover is given by the curve $y^6=x(x-1)(x-t)$ (compare \cite{FM08}).
\par
If either of $a$ or $b$ is different from $1/2$, we may suppose
that $b \neq 1/2$, maybe after flipping along the first diagonal.
In this case $\alpha \beta \gamma_1^{-1} \alpha$ also represents
a core curve of a cylinder of slope $1/2$. Since the $\varphi_\opt$-image
of this element is trivial, we obtain the desired contradiction. 
\hfill $\Box$
\par
With Corollary~\ref{cor:nocase1} and Lemma~\ref{le:case2} 
we have completed the proof of Theorem~\ref{Main}. 
\par
In genus $5$ the only possible location of a ST-curves are the finitely
many square-tiled coverings given by the cases (3), (4), (5), (7), (8) and (9)
and possibly infinitely many square-tiled coverings in the case (5). However
Corollary~\ref{eqlenalldir} implies strong restrictions
to the possible monodromy representations of such a ST-curve in $g=5$ -- if it exists.
\par

\section{Lyapunov exponents of Teichm\"uller curves} \label{sec:Lyap}

Fix an $\SL_2(\RR)$-invariant, ergodic measure $\mu$ on $\Omega M_g^*$.
The Lyapunov exponents for the Teichm\"uller geodesic flow
on $\Omega M_g$ measure the logarithm of the growth rate of the 
Hodge norm of cohomology classes during parallel transport along 
the geodesic flow. More precisely, let $\VV$ be the restriction
of the real Hodge bundle $V$ (i.e.\ the bundle with fibers $H^1(X,\RR)$)
to the support of $\mu$. Let $S_t$ be the lift of the geodesic
flow to $V$ via the Gauss-Manin connection. Then Oseledec's theorem 
shows the existence of a filtration
$$V = V_{\lambda_1} \supset \cdots V_{\lambda_k} \supset 0$$
by measurable vector subbundles such that, for almost all $m \in M$ 
and all $v \in V_m \sms \{0\}$, one has
$$||S_t(v)|| = {\rm exp}(\lambda_i t + o(t)),$$
where $i$ is the maximal value such that
$v\in (V_i)_m$.
The numbers $\lambda_i$ for $i=1,\ldots,k\leq {\rm rank}(V)$
are called  the {\em  Lyapunov exponents of $S_t$}. Note
that these exponents are unchanged if we replace
the support of $\mu$ by a finite unramified covering with
a lift of the flow and the pullback of $V$. We adopt the convention
to repeat the exponents according to the
rank of $V_i/V_{i+1}$ such that we will always
have $2g$ of them, possibly some of them equal.
Since $V$ is symplectic, the spectrum is symmetric, i.e.\
$\lambda_{g+k} = - \lambda_{g-k+1}$. 
The reader may consult \cite{Fo06} or \cite{Zo06} for a more 
detailed introduction to the subject. The guiding questions
are existence of zero Lyapunov exponents and simplicity
of the set of Lyapunov exponents.
\par
While the Lyapunov spectrum defined is for any $\SL_2(\RR)$-invariant
measure $\mu$, the extremal
case that the support of $\mu$ is very large, namely a stratum,
is best understood. It has been shown by Veech (\cite{Ve86}) that 
the second Lyapunov exponent is strictly smaller than one. Forni has
shown in \cite{Fo02}, that $\lambda_g >0$. Subsequently,
Zorich's conjecture, saying that the Lyapunov spectrum is simple, was 
completed by \cite{AV07}.
\par
The analogous statements fail when the support of $\mu$ is minimal, i.e.\ when 
$\mu$ is supported on
a closed orbit or equivalently on the natural lift of a Teichm\"uller curve
to $\Omega M_g$. In this case it was discovered by Forni (\cite{Fo06})
that there exists an example with zero Lyapunov
exponent and even an example where the second Lyapunov exponent $\lambda_2$
is zero, see also below. In this case we say that the Lyapunov 
spectrum is {\em totally degenerate}.
\par
We start we a summary of known results on the Lyapunov spectrum for 
Teichm\"uller curves.
\par
\begin{Thm} [\cite{BM05}]
Let $f: X \to C$ be the universal family over a Teichm\"uller curve 
with the conventions stated above. Suppose the VHS has a rank two 
summand $\LL$ whose $(1,0)$-part of the Hodge 
filtration is given by a line bundle $\cLL$, then the Lyapunov spectrum of the 
Teichm\"uller curve is the union of plus and minus
$$\lambda = \deg \cLL/(g-2+ \#\{\ol{C} \setminus C\}/2)$$
and the spectrum of the complement of $\LL$.
\end{Thm}
\par
Sometimes e.g.\ for cyclic coverings of the projective line, 
this result allows to calculate all the Lyapunov exponents. See
\cite{BM05} Proposition~3.4 for the precise statement. Note that
loc.\ cit.\ and the above theorem applies only to the summands of the VHS
that are non-unitary. For summands of the VHS that are unitary, e.g.\ those 
for which the Hodge filtration is trivial, the Lyapunov exponents are zero.
\par
We recall that a family of abelian varieties $g: \cYY  \to C$ 
has a fixed part of dimension $r$, if there is an abelian variety
$A$ of dimension $r$ and an injection $A \times C \to \cYY$
over the identity of $C$.
\par
\begin{Prop}
Suppose that $f$ has a fixed part of dimension $r$. Then $2r$ of the
Lyapunov exponents are zero. 
\end{Prop}
\par
\begin{proof}
Under this assumptions, the variation of Hodge structures
splits off a trivial local system $\UU$ of rank $2r$. All the
Lyapunov exponents corresponding to $\UU$ are obviously zero.
\end{proof}
\par
This together with the results of Section~\ref{sec:cyccov} imply:
\par
\begin{Cor}[\cite{Fo06}, \cite{FM08}] \label{cor:FMtotdeg}
The Lyapunov spectrum of the family of curves
$$y^4 = x (x-1)(x-t) $$
resp.\
$$y^6 = x (x-1)(x-t) $$
is totally degenerate, that is, it equals
$$ \{\pm 1,0,0\} \quad \text{resp.} \quad \{ \pm 1,0,0,0 \}.$$ 
\end{Cor}
\par
\begin{Prop}
Suppose that the Lyapunov spectrum of a Teichm\"uller curve $f$ is 
totally degenerate. Then $f$ is a ST-curve.
\end{Prop}
\par
\begin{proof}
We write the VHS of $f$ as $\VV = \LL \oplus \MM$, where
$\LL$ is the maximal Higgs subsystem, i.e.\ contributes the Lyapunov 
exponents $\pm 1$ and $\MM$ its orthogonal complement. The hypothesis
implies that all the Lyapunov exponents of $\MM$ are zero. We claim, 
that consequently, there is no expansion of the Hodge norm on $\MM$
along the geodesic flow. Since the geodesic flow acts ergodically
on the unit tangent bundle to $C$, the claim implies that 
the local system $\MM$ is unitary.
Since Teichm\"uller curves are non-compact, one shows as in 
\cite{ViZu04} that the splitting of $\VV$ is defined over $\QQ$, i.e.\
a splitting of abelian varieties up to isogeny. Since the local 
system $\VV$ carries a $\ZZ$-structure, $\MM$
is local system with transition maps in $\Sp(2g-2,\ZZ[1/N])$ for
some $N \in \NN$. The intersection of $\Sp(2g-2,\ZZ[1/N])$ and
the unitary group is obviously finite. Thus, after a finite
base change, the family $f$ has a fixed part of dimension $g-1$.
Consequently, it is a ST-curve.
\par
To prove the claim, we have to show that the vanishing of the
Lyapunov exponents -- which a priori implies just subexponential
growth -- implies that the Hodge norms remain indeed constant.
For this purpose we use \cite{Fo02} Corollary 5.3. It states that
that
$$ \lambda_2 + \cdots \lambda_g = \int_{{\rm supp}(\mu)} \Lambda_2(q) + \cdots
\Lambda_g(q) \mu(q),$$
where the $\Lambda_i(q) \geq 0$ are the eigenvector of a hermitian
form $H_q$ obtained as a Hodge norm modified by $q$. If $\lambda_i=0$
for all $i>1$ then $\Lambda_i = 0$ $\mu$-almost everywhere and
for all $i>1$. By  \cite{Fo02} Lemma~2.1, for all cohomology classes $c$
in the orthogonal complement of $\langle \omega, \bar{\omega} \rangle$,
the derivative of the Hodge norm is given by
$$ \frac{d}{dt}||S_t(c)|| = -2 {\rm Re}(B_q(S_t(c), S_t(c))),$$
where $B_q$ is a bilinear form with $H_q = B_q^* B_q$. Since
$H_q$ vanishes, so does $B_q$ and the derivative.
\end{proof} 
\par
We may thus restate the main theorem as follows.
\begin{Cor}
Up to possible exceptions in $M_5$ the only Teichm\"uller curves
with totally degenerate  Lyapunov spectrum are the
ones listed in Corollary~\ref{cor:FMtotdeg}.
\end{Cor}
\par

\par
Martin M{\"o}ller: Max-Planck Institut f\"ur Mathematik\newline
Postfach 7280 \newline
53072 Bonn, Germany \newline
e-mail: moeller@mpim-bonn.mpg.de 
\par

\end{document}